\documentclass[11pt,reqno]{amsart}

\usepackage{amsmath, amsfonts, amsthm, amssymb, color, cite}
\newtheorem{theorem}{Theorem}[section]
\newtheorem{lemma}{Lemma}[section]

\theoremstyle{definition}
\newtheorem{definition}{Definition}[section]

\theoremstyle{remark}
\newtheorem{remark}{Remark}[section]

\numberwithin{equation}{section}
\allowdisplaybreaks

\newcommand{\R}{{\mathbb R}}

\def\f{\frac}

\def\hf1{^\f{1}{1-\xi^2}}

\def\be{\begin{equation}}
\def\en{\end{equation}}
\def\bs{\begin{split}}
\def\es{\end{split}}
\def\ba{\begin{align}}
\def\ea{\end{align}}

\newcommand{\eps}{\varepsilon}
\author[Wentao Cao]{Wentao Cao}
\address{Institute f\"{u}r mathematik, Universit\"{a}t Leipzig, D-04109, Leipzig, Germany}
\email{wentao.cao@math.uni-leipzig.de}
\author[Feimin Huang]{Feimin Huang}
\address {School of Mathematical Sciences, University of Chinese Academy of Sciences, Beijing 100049, China; Academy of Mathematics and System Sciences, Chinese Academy of Sciences, Beijing 100190, China}
\email{fhuang@amt.ac.cn}
\author[Difan Yuan]{Difan Yuan}
\address{University of Chinese Academy of Sciences, Institute of Applied Mathematics,
AMSS, Beijing 100190, China
}
\email{yuandf@amss.ac.cn}

\title[General Nozzle Flow]
{Global Entropy Solutions to the Gas Flow in General Nozzle}


\date{\today}

\begin{document}

\begin{abstract}
We are concerned with the global existence of entropy  solutions for the compressible Euler equations describing the gas flow in a nozzle with general cross-sectional area, for both isentropic and isothermal fluids. New viscosities are delicately designed to obtain the uniform bound of approximate solutions. The vanishing viscosity method and compensated compactness framework are used to prove the convergence of approximate solutions. Moreover, the entropy solutions for both cases are uniformly bounded independent of time. No smallness condition is assumed on initial data. The techniques developed here can be applied to compressible Euler equations with general source terms.
\end{abstract}
\maketitle
\medskip

\noindent {\bf 2010 AMS Classification}:  35L45, 35L60, 35Q35.

\section{Introduction}
%

We consider one dimensional gas flow in a general nozzle for the isentropic and isothermal flows separately. The nozzle is widely used in some types of steam turbines, rocket engine nozzles, supersonic jet engines, and jet streams in astrophysics. The motion of the nozzle flow is governed by the following system of compressible Euler equations:
\begin{eqnarray}\label{iso1}
\left\{ \begin{array}{ll}
\displaystyle \rho_t+m_x=a(x)m,\,x\in\R,\,t>0,\\
\displaystyle  m_t+\left(\frac{m^2}{\rho}+p(\rho)\right)_x=a(x)\frac{m^2}{\rho},
\,x\in\R,\,t>0,
\end{array}
\right.
\end{eqnarray}
where $\rho$ is the density, the momentum $m=\rho u$ with $u$ being the velocity, and $p(\rho)$ is the pressure of the gas. Here the given function $a(x)$ is represented by $a(x)=-\frac{A'(x)}{A(x)}$ with  $A(x)\in C^2(\R)$ being a slowly variable cross-sectional area at $x$ in the nozzle.  For $\gamma$-law gas, $p(\rho)=p_0\rho^\gamma$ with  $\gamma$ denoting the adiabatic exponent and $p_0=\frac{\theta^{2}}{\gamma},\theta=\frac{\gamma-1}{2}$. When $\gamma>1,$ \eqref{iso1} is called the isentropic gas flow. When $\gamma=1,$  \eqref{iso1} is called isothermal one.
We consider the Cauchy problem for \eqref{iso1} with large initial data
\begin{equation}\label{ini1}
(\rho, m)|_{t=0}=(\rho_0(x), m_0(x))\in L^\infty.
\end{equation}
The above Cauchy problem \eqref{iso1}-\eqref{ini1} can be written in compact form as follows:
\begin{eqnarray}\label{iso3}
\left\{ \begin{array}{llll}
\displaystyle U_t+f(U)_x=g(x, U),\\
\displaystyle U|_{t=0}=U_0(x), x\in\R,
\end{array}
\right.
\end{eqnarray}
where $U=(\rho, m)^\top,$ $f(U)=(m, \frac{m^2}{\rho}+p(\rho))^\top,$ and
$g(x, U)=(-\frac{A'(x)}{A(x)}m, -\frac{A'(x)}{A(x)}\frac{m^2}{\rho})^\top.$

There have been extensive studies and applications of homogeneous $\gamma$-law gas, i.e., $g(x,U)=0$. Diperna \cite{Diperna} proved the global existence of entropy solutions with large initial data by the theory of compensated compactness and vanishing viscosity method for $\gamma=1+\frac{2}{2n+1},$ where $n$ is a positive integer. Subsequently, Ding, Chen, and Luo\cite{Ding,DingDing} and Chen \cite{chen} successfully extended the result to $\gamma\in(1, \frac{5}{3}]$ by using a Lax-Friedrichs scheme. Lions, Perthame, and Tadmor \cite{Lions1} and Lions, Perthame, and Souganidis \cite{Lions2} treated the case $\gamma>\frac{5}{3}$. The existence of entropy solutions to the isothermal gas, i.e., $\gamma=1$,  was proved in Huang and Wang \cite{HuangWang} by introducing complex entropies and utilizing  the analytic extension method.

For the isentropic Euler equations with source term, Ding, Chen, and Luo \cite{Ding1} established a general framework to investigate the global existence of entropy solution through the fractional step Lax-Friedrichs scheme and compensated compactness method. Later on, there have been extensive studies on the inhomogeneous case (see \cite{Chen2,Chen3,Lu1997, Marcati, Marcati2, Wangzejun, Wangzejun2}).
For the nozzle flow problem, see \cite{Courant,Embid,Glaz,Glimm1984,Liu1979,Liu1982,Liu1987,Whitham}. For converging-diverging de Laval nozzles, as flow speed accelerates from the subsonic to the supersonic regime, the physical properties of nozzle and diffuser flows are altered. This kind of nozzle is particularly designed to converge to a minimum cross-sectional area and then expand. Liu \cite{Liu1979} first proved the existence of a global solution with initial data of small total variation and away from sonic state by a Glimm scheme. Tsuge \cite{Tsuge3,Tsuge4,Tsuge5} first studied the global existence of solutions for Laval nozzle flow and transonic flow for large initial data by introducing a modified Godunov scheme. Recently, Chen and Schrecker \cite{Chen2018} proved the existence of globally defined entropy solutions in transonic nozzles in an $L^p$ compactness framework, whose uniform bound of approximate solutions may depend on time $t$. In our paper, we are focusing on the $L^{\infty}$ compactness framework. Moreover, general cross-sectional areas of nozzles are considered, which include several important physical models, such as the de Laval nozzles with closed ends, that is, the cross-sectional areas are tending to zero as $x\rightarrow\infty.$
%

In our paper, we assume the cross-sectional area function $A(x)>0$ satisfies that there exists a $C^{1, 1}$ function $a_0(x)\in L^1(\R)$ such that
\begin{equation}\label{acondition}
\left|\frac{A'(x)}{A(x)}\right|=|a(x)|\leq a_0(x).
\end{equation}
Here, $A(x)>0$ is a natural assumption. The smallest cross-sectional area of the nozzle is the throat of the nozzle. We allow the general varied cross-sectional area and no assumption is assumed on the sign of $a(x).$

The main purpose of the present paper is to prove the existence of a global entropy solution with uniform bound independent of time for large initial data in both the isentropic case $1<\gamma<3$ and isothermal one $\gamma=1$. We are interested in solutions that can reach the vacuum $\rho=0.$ Near the vacuum, the system \eqref{iso1}-\eqref{ini1} is degenerate and  the velocity $u$ cannot be defined uniquely. We define the weak entropy solution as follows.
\begin{definition}
A measurable function $U(x, t)$ is called a global weak solution of the Cauchy problem \eqref{iso3} if
\begin{equation*}
\int_{t>0}\int_\R U\varphi_t+f(U)\varphi_x+g(x, U)\varphi dxdt+\int_\R U_0(x)\varphi(x, 0)dx=0
\end{equation*}
holds for any test function $\varphi\in C^1_0(\R\times\R^+)$. In addition,
for the isentropic flow,  if $U$ also satisfies that for any weak entropy pair $(\eta, q)$ (see Section 2), the inequality
\begin{equation}\label{entropyinequal}
\eta(U)_t+q(U)_x-\nabla\eta(U)\cdot g(x,U)\leq0
\end{equation}
holds in the sense of distributions, then $U$ is called a weak entropy solution to \eqref{iso3}. For the isothermal flow,  $U$ is called a weak entropy solution if $U$ additionally satisfies \eqref{entropyinequal} for mechanical entropy pair
$$\eta_*=\frac{m^2}{2\rho}+\rho\ln\rho,\quad  q_*=\frac{m^3}{2\rho^2}+m\ln\rho.$$
\end{definition}

Two main results of the present paper are given as follows.
\begin{theorem}\label{mainisen}{\text(isentropic case)}
Let $1<\gamma<3.$ Assume that there is a positive constant $M$ such that the initial data satisfies
\begin{equation*}
0\leq\rho_0(x)\leq M,~~|m_0(x)|\leq M\rho_0(x),~\text{ $a.e.,x\in\R $},
\end{equation*}
and $a(x)$ satisfies \eqref{acondition} with
\begin{equation}\label{a01}
\|a_0(x)\|_{L^{1}(\R)}\leq\frac{1-\theta}{1+\theta}.
\end{equation}
Then, there exists a global entropy solution of \eqref{iso1}-\eqref{ini1} satisfying
\begin{equation*}
0\leq\rho(x,t)\leq C,~|m(x, t)|\leq C\rho(x, t), \text{ a.e., $(x,t)\in\R\times\R^{+}$},
\end{equation*}
where $C$ depends only on initial data and is independent of time $t$.
\end{theorem}
\begin{theorem}\label{main}{\text(isothermal case)}
Let $\gamma=1.$ Assume that there is a positive constant $M$ such that the initial data satisfy
\begin{equation*}
0\leq\rho_0(x)\leq M,~~|m_0(x)|\leq\rho_0(x)(M+|\ln\rho_0(x)|),~\text{ a.e., }x\in\R ,
\end{equation*}
and  $a(x)$ satisfies \eqref{acondition} with
\begin{equation}\label{a02}
\|a_0(x)\|_{L^{1}(\R)}\leq\tfrac{1}{2}.
\end{equation}
Then, there exists a global entropy solution of \eqref{iso1}-\eqref{ini1} satisfying
\begin{equation*}
0\leq\rho(x,t)\leq C,~|m(x, t)|\leq\rho(x, t)(C+|\ln\rho(x, t)|), \text{ a.e., }x\in\R,
\end{equation*}
where $C$ depends only on initial data and is independent of time $t$.
\end{theorem}%
\begin{remark}
Here, the conditions \eqref{a01} (Theorem \ref{mainisen}) and \eqref{a02} (Theorem \ref{main}) are assumed to guarantee a uniform bound of $(\rho, m)$ independent of time. This condition illustrates a new physical phenomena that is important in engineering. For example, if we consider an isothermal nozzle with a monotone cross-sectional area, $a_0(x)=\frac{A'(x)}{A(x)}\geq0,$ and denote $A_+$ and $A_-$ the far field of a variable cross-sectional area, respectively, then the ratio of the outlet and inlet cross-sectional area can be controlled, i.e., $\frac{A_+}{A_-}\leq e^{\frac{1}{2}}.$
\end{remark}
\begin{remark}
The condition \eqref{a01} in Theorem \ref{mainisen} is different from that in Tsuge \cite{Tsuge5}. Here, in our paper, we allow $1<\gamma<3.$
\end{remark}
 The main difficulty we came across is how to construct approximate solutions with uniform bound independent of time. Another difficulty is the interaction of nonlinear resonance between the characteristic modes and geometrical source terms. Our strategy is applying the maximum principle (Lemma 3.1) introduced in \cite{Huang20171,Huang20172}, which is similar to invariant region theory \cite{Smoller}, to a viscous equation with novel viscosity. To be more specific, for the isentropic case, we add $-2\varepsilon b(x)\rho_x$ on the momentum equation (c.f \eqref{isen-vis}); for the isothermal case, we raise $n:=A(x)\rho$ with $\delta$ and also add $-4\varepsilon b(x)n_x$ on the momentum equation (c.f \eqref{iso-vis2}). Two modified Riemann invariants are introduced and a system of decoupled new parabolic equations along the characteristic are derived. Owing to the hyperbolicty structure of \eqref{iso1}, we can transform the integral of source terms along characteristics with time $t$ into the integral with space $x.$ Finally after establishing the estimate of $H^{-1}_{loc}$ compactness, we apply a compensated compactness framework in \cite{Ding,DingDing, Lions2, HuangWang} to show the convergence of approximate solutions. To the best of our knowledge, for the isothermal flow, the uniform bound for the approximate solutions depends on time $t$ in all the previous results. We remark that the method in our paper can be applied to obtain the existence of weak solutions of related gas dynamic models, such as Euler-Poisson for a semiconductor model \cite{Huang20172} or an Euler equation with geometric source terms \cite{Huang20171}, and may also shed light on the large time behavior of entropy solutions. Besides, we avoid a laborious numerical scheme to construct approximate solutions.

The present paper is organized as follows. In Section \ref{formula}, we introduce some basic notions and formulas for the isentropic Euler system. In Section \ref{isentheorem}, we prove Theorem \ref{mainisen} for the global existence of isentropic gas flow in general nozzle.  Subsequently, in Section \ref{formula2}, we further formulate several preliminaries and formula for the isothermal Euler system. The proof of Theorem \ref{main} for global existence of isothermal gas flow in general nozzle will be presented in Section \ref{isotheorem}. In the appendix, we provide the proof of variant version of invariant region theory for completeness.
\section{Preliminary and Formulation for Isentropic Flow}\label{formula}
First we list some basic notation for the isentropic system \eqref{iso1}.
The eigenvalues are
\begin{equation*}
\lambda_1=\frac{m}{\rho}-\theta\rho^\theta,\quad
\lambda_2=\frac{m}{\rho}+\theta\rho^\theta,\quad
\end{equation*}
and the corresponding right eigenvectors are
\begin{equation*}
r_1=\left[\begin{array}{cc}
1\\ \lambda_1
\end{array}
\right],\quad
r_2=\left[\begin{array}{cc}
1\\ \lambda_2
\end{array}
\right].
\end{equation*}
The Riemann invariants $w, z$ are given by
\begin{equation}\label{2.3}
w=\frac{m}{\rho}+\rho^\theta,\quad z=\frac{m}{\rho}-\rho^\theta,
\end{equation}
satisfying $\nabla w\cdot r_1=0$ and $\nabla z\cdot r_2=0$.  A pair of functions $(\eta, q):\R^+\times\R\mapsto\R^2$ is defined to be an entropy-entropy flux pair if it satisfies
\begin{equation*}
\nabla q(U)=\nabla\eta(U)\nabla\left[\begin{array}{ccc}
m\\ \frac{m^2}{\rho}+p(\rho)
\end{array}
\right].
\end{equation*}
When $$\eta\left|_{\frac{m}{\rho}\text{ fixed }}\rightarrow0, \text{ as } \rho\rightarrow 0, \right.$$ $\eta(\rho, m)$ is called weak entropy.
In particular, the mechanical entropy pair
\begin{equation*}
\eta^*(\rho, m)=\frac{m^2}{2\rho}+\frac{p_0\rho^\gamma}{\gamma-1},~~
q^*(\rho, m)=\frac{m^3}{2\rho^2}+\frac{\gamma p_0\rho^{\gamma-1}m}{\gamma-1}
\end{equation*}
is a strictly convex entropy pair.
As shown in \cite{Lions1} and \cite{Lions2}, any weak entropy for the system \eqref{iso1} is given by
\begin{equation}\label{2.6}
\begin{split}
\eta=\rho\int_{-1}^1\chi(\frac{m}{\rho}+\rho^\theta s)(1-s^2)^\lambda ds,~~
q=\rho\int_{-1}^1(\frac{m}{\rho}+\rho^\theta\theta s)\chi(\frac{m}{\rho}+\rho^\theta s)(1-s^2)^\lambda ds
\end{split}
\end{equation}
with  $\lambda=\frac{3-\gamma}{2(\gamma-1)}$ for any function  $\chi(\cdot)\in C^2(\R)$.

\section{Proof of Theorem \ref{mainisen}}\label{isentheorem}
\subsection{Construction of approximate solutions}
We first construct approximate solutions to \eqref{iso1} satisfying the framework in \cite{Ding,DingDing, Lions2}. Indeed,  for any $\eps\in(0,1)$ we construct approximate solutions by adding suitable artificial viscosity as follows:
\begin{eqnarray}\label{isen-vis}
\left\{ \begin{array}{llll}
\displaystyle \rho_t+m_x=a(x)m+\eps\rho_{xx},\\
\displaystyle m_t+\left(\frac{m^2}{\rho}+p(\rho)\right)_x
=a(x)\frac{m^2}{\rho}+\eps m_{xx}-2\eps b(x)\rho_x
\end{array}
\right.
\end{eqnarray}
with initial data
\begin{equation}\label{isenini-vis}
(\rho, m)|_{t=0}=(\rho_0^\eps(x), m_0^\eps(x))=(\rho_0(x)+\eps, m_0(x))\ast j^\eps,
\end{equation}
where $b(x)$ is a function to be given later, and $j^\eps$ is the standard mollifier.

\subsection{Global existence of approximate solutions}\label{approximate}

For the global existence to Cauchy problem \eqref{isen-vis}-\eqref{isenini-vis}, we have the following.
\begin{theorem}\label{thm-isenthvis}
For any time $T>0,$ there exists a unique global classical bounded solution to the Cauchy problem \eqref{isen-vis}-\eqref{isenini-vis} that has following $L^{\infty}$ estimates
\begin{equation}\label{bound}
e^{-C(\eps, T)}\leq \rho^\eps(x, t)\leq C, ~~|m^\eps(x, t)|\leq C\rho^\eps(x, t).
\end{equation}
\end{theorem}
 We shall show Theorem \ref{thm-isenthvis} in two steps.  In the section, we omit the upper index $\eps$ for simplicity.

\textbf{Step 1. Uniform upper bound.}  First, we can rewrite the first equation of \eqref{isen-vis} as
\begin{equation*}
\rho_t+u \rho_x=\eps\rho_{xx}+\rho(a(x)u-u_x),
\end{equation*}
and then applying the maximum principle of parabolic equation yields that
\begin{align*}
\rho\geq\min\rho_0(x)e^{-\int_0^t\|a(x)u-u_x\|_{L^\infty}ds}>0,
\end{align*}
which implies $w\geq z.$   Second, we recall a revised version of the invariant region theory \cite{Smoller} introduced in \cite{Huang20171,Huang20172}.
\begin{lemma}(Maximum principle)\label{modified maximum}
Let $p(x,t), q(x,t)$, $(x,t)\in  \R\times[0,T]$ be any bounded classical solutions of the quasilinear parabolic system
\begin{eqnarray}\label{pq}
\left\{ \begin{aligned}
\displaystyle &p_t+\mu_1 p_x=
\eps p_{xx}+a_{11}p+a_{12}q+R_1,\\
\displaystyle &q_t+\mu_2 q_x=
\eps q_{xx}+a_{21}p+a_{22}q+R_2
\end{aligned}
\right.
\end{eqnarray}
with initial data
$p(x,0)\leq0, ~~q(x, 0)\geq0, $
where $$\mu_{i}=\mu_i(x,t,p(x,t),q(x,t)),a_{ij}=a_{ij}(x,t,p(x,t),q(x,t)),$$ and the source terms $$R_i=R_i(x,t,p(x,t),q(x,t),p_{x}(x,t),q_{x}(x,t)),i,j=1,2,\forall(x,t)\in\R\times[0,T]$$ $\mu_{i},a_{ij}$are bounded with respect to $(x,t,p,q)\in\R\times[0,T]\times K,$ where $K$ is an arbitrary compact subset in $\R^2,$  $a_{12},a_{21},R_{1},R_{2}$ are continuously differentiable with respect to $p,q.$
Assume the following conditions hold:
\begin{description}
  \item[(C1)]  When $p=0$ and $q\geq0,$ there is $a_{12}\leq0;$ when $q=0$ and $p\leq0,$ there is   $a_{21}\leq0.$
  \item[(C2)]  When $p=0$ and $q\geq0,$ there is $~R_1=R_1(x,t,0,q,\zeta,\eta)\leq0;$ when $q=0$ and $p\leq0,$ there is $~R_2=R_2(x,t,p,0,\zeta,\eta)\geq0.$
\end{description}
Then for any $(x, t)\in\R\times[0,T],$
$p(x,t)\leq0, ~~q(x, t)\geq0.$
\end{lemma}
\begin{remark}
The modified version of invariant region theory (Lemma \ref{modified maximum}) is valid not only for the Cauchy problem with source terms, but also for the initial boundary value problem with Dirichlet and Neumann boundary conditions.
\end{remark}
We shall apply maximum principle Lemma \ref{modified maximum} to get the uniform bound of $\rho, m$. By the formulas of Riemann invariants \eqref{2.3}, the viscous perturbation system \eqref{isen-vis} can be transformed as
\begin{eqnarray}\label{wzisen}
\left\{ \begin{aligned}
\displaystyle &w_t+\lambda_2 w_x=\eps w_{xx}+2\eps(w_x-b)\frac{\rho_x}{\rho}-\eps\theta(\theta+1)\rho^{\theta-2}\rho_x^2+\theta\frac{w^2-z^2}{4}a(x),\\
\displaystyle &z_t+\lambda_1 z_x=\eps z_{xx}+2\eps(z_x-b)\frac{\rho_x}{\rho}+\eps\theta(\theta+1)\rho^{\theta-2}\rho_x^2-\theta\frac{w^2-z^2}{4}a(x).
\end{aligned}
\right.
\end{eqnarray}
Set the control functions $(\phi,\psi)$ as
\begin{equation*}
\begin{split}
&\phi=C_0+\eps\|b'(x)\|_{L^{\infty}}t+\int_{-\infty}^xb(y)dy,\\
&\psi=C_0+\eps\|b'(x)\|_{L^{\infty}}t+\int^{\infty}_xb(y)dy.
\end{split}
\end{equation*}
Then a simple calculation shows that
\begin{equation*}
\begin{split}
&\phi_t=\eps\|b'(x)\|_{L^{\infty}},~\phi_x=b(x), ~\phi_{xx}=b'(x);\\
&\psi_t=\eps\|b'(x)\|_{L^{\infty}},~\psi_x=-b(x), ~\psi_{xx}=-b'(x).
\end{split}
\end{equation*}
Define the modified Riemann invariants $(\bar{w},\bar{z})$ as
\begin{equation}\label{r}
\bar{w}=w-\phi, ~~\bar{z}=z+\psi.
\end{equation}
Inserting \eqref{r} into \eqref{wzisen} yields the decoupled equations for $\bar{w}$ and $\bar{z}:$
\begin{eqnarray}\label{phipsi11}
\left\{\begin{aligned}
\bar{w}_t+\lambda_2\bar{w}_x=&\eps\bar{w}_{xx}
+\eps\phi_{xx}-\phi_t-\lambda_2\phi_x+2\eps\frac{\rho_x}{\rho}\bar{w}_x-\eps\theta(\theta+1)\rho^{\theta-2}\rho_x^2 \\
&+\theta\frac{(\bar{w}+\phi)^{2}-(\bar{z}-\psi)^{2}}{4}a(x),\\
\bar{z}_t+\lambda_1\bar{z}_x=&\eps\bar{z}_{xx}
-\eps\psi_{xx}+\psi_t+\lambda_1\psi_x+2\eps\frac{\rho_x}{\rho}\bar{z}_x+\eps\theta(\theta+1)\rho^{\theta-2}\rho_x^2\\
&-\theta\frac{(\bar{w}+\phi)^{2}-(\bar{z}-\psi)^{2}}{4}a(x).
\end{aligned}
\right.
\end{eqnarray}
Noting that
\begin{equation*}
\begin{split}
&\lambda_1=\frac{w+z}{2}-\theta\frac{w-z}{2},\\
&\lambda_2=\frac{w+z}{2}+\theta\frac{w-z}{2},
\end{split}
\end{equation*}
the system \eqref{phipsi11} becomes
\begin{eqnarray}\label{rst}
\displaystyle\left\{ \begin{aligned} &\bar{w}_t+(\lambda_2-2\eps\frac{\rho_x}{\rho})\bar{w}_x
=\eps\bar{w}_{xx}+a_{11}\bar{w}
+a_{12}\bar{z}+R_1,\\
&\bar{z}_t+(\lambda_1-2\eps\frac{\rho_x}{\rho})\bar{z}_x
=\eps\bar{z}_{xx}+a_{21}\bar{w}
+a_{22}\bar{z}+R_2,
\end{aligned}
\right.
\end{eqnarray}
where
\begin{equation*}
\begin{split}
&a_{11}=-\left(\frac{1+\theta}{2}\phi_x-\theta\frac{\bar{w}+2\phi}{4}a(x)\right), \quad
a_{12}=-\left(\frac{1-\theta}{2}\phi_x+\theta\frac{\bar{z}-2\psi}{4}a(x)\right),\\
&a_{21}=\left(\frac{1-\theta}{2}\psi_x-\theta\frac{\bar{w}+2\phi}{4}a(x)\right), \quad
a_{22}=\left(\frac{1+\theta}{2}\psi_x+\theta\frac{\bar{z}-2\psi}{4}a(x)\right),
\end{split}
\end{equation*}
and
\begin{equation*}
\begin{split}
&R_1=\eps\phi_{xx}-\phi_t-\frac{1+\theta}{2}\phi\phi_{x}+\frac{1-\theta}{2}\psi\phi_{x}
-\eps\theta(\theta+1)\rho^{\theta-2}\rho_x^2+\theta\frac{\phi^{2}-\psi^{2}}{4}a(x),\\
&R_2=-\eps\psi_{xx}+\psi_t+\frac{1-\theta}{2}\phi\psi_{x}-\frac{1+\theta}{2}\psi\psi_{x}
+\eps\theta(\theta+1)\rho^{\theta-2}\rho_x^2-\theta\frac{\phi^{2}-\psi^{2}}{4}a(x).
\end{split}
\end{equation*}
To apply Lemma \ref{modified maximum}, we need to verify $\textbf{(C1)}$ and $\textbf{(C2)}.$   For $\textbf{(C1)},$
when $\bar{w}=0 ,\bar{z}\geq0,$ we have
\begin{align*}
0\leq\bar{z}=z+\psi\leq w+\psi=\phi+\psi,
\end{align*}
and then
\[ \begin{aligned}
 a_{12}=&-\frac{1-\theta}{2}\left(b(x)+\frac{\theta}{2(1-\theta)}(\bar{z}-2\psi)a(x)\right) \\
 \leq& \displaystyle\left\{
 \begin{aligned}
&-\frac{1-\theta}{2}\left(b(x)-\frac{\theta}{2(1-\theta)}(\phi-\psi)|a(x)|\right) \text{ if} ~~a(x)<0,\\
&-\frac{1-\theta}{2}\left (b(x)-\frac{\theta}{2(1-\theta)}2\psi|a(x)|\right)  \text{ if} ~~a(x)\geq0.
\end{aligned}
\right.
\end{aligned} \]
Hence, we take $b(x)=M_0a_0(x)$ with
$$M_0\geq\frac{\theta}{2(1-\theta)}\max(\phi-\psi,2\psi),$$
and using \eqref{acondition}, one has $a_{12}\leq0.$ Moreover, when $\bar{w}\leq0, \bar{z}=0,$ we have
\begin{align*}
0\geq\bar{w}=w-\phi\geq z-\phi=-\psi-\phi,
\end{align*}
and then
\[ \begin{aligned}
 a_{21}=& -\frac{1-\theta}{2}\left(b(x)+\frac{\theta}{2(1-\theta)}(\bar{w}+2\phi)a(x)\right) \\
\leq &\displaystyle\left\{ \begin{aligned}
& -\frac{1-\theta}{2}\left(b(x)-\frac{\theta}{2(1-\theta)}2\phi|a(x)|\right) \text{ if} ~~a(x)<0,\\
& -\frac{1-\theta}{2}\left(b(x)-\frac{\theta}{2(1-\theta)}(\psi-\phi)|a(x)|\right) \text{ if} ~~a(x)\geq0,
\end{aligned}
\right. \\
\leq &~0,
\end{aligned} \]
provided $$M_0\geq\frac{\theta}{2(1-\theta)}\max(\phi-\psi,2\phi).$$
Thus we require\\
 $$M_0\geq\frac{\theta}{2(1-\theta)}\max
 \left(\int_{-\infty}^xb(y)dy-\int_x^\infty b(y)dy,2C_0+2\eps\|b'(x)\|_{L^{\infty}}t+2\int_{-\infty}^\infty b(y)dy\right).$$
Taking $\eps$ sufficient small such that $\eps\|b'(x)\|_{L^{\infty}}T\leq1,$ we have
$$M_0\geq\frac{\theta}{1-\theta}\left({C_0+1+M_0\|a_0(x)\|_{L^1}}\right),$$
that is,
\begin{equation}\label{2}
\|a_0(x)\|_{L^1}\leq\frac{1-\theta}{\theta}-\frac{C_0+1}{M_0}.
\end{equation}
Hence $\textbf{(C1)}$ is also satisfied by $(\bar{w}, \bar{z}).$  As for $\textbf{(C2)},$ one can derive
\begin{equation*}
\begin{split}
R_1&\leq \eps b'(x)-\eps\|b'(x)\|_{L^{\infty}}\\
&\qquad +b(x)\left(-\frac{1+\theta}{2}\phi+\frac{1-\theta}{2}\psi\right)+\theta\frac{(\phi+\psi)(\phi-\psi)}{4}a(x)\\
&\leq b(x)\left(-\theta C_0-\eps\theta\|b'(x)\|_{L^{\infty}}t-\frac{1+\theta}{2}\int_{-\infty}^xb(y)dy
+\frac{1-\theta}{2}\int^{\infty}_xb(y)dy\right)\\
&\qquad +\frac{\theta}{4}({2C_0+2\eps\|b'(x)\|_{L^{\infty}}}t+\|b\|_{L^1})
\left(\int_{-\infty}^xb(y)dy-\int_x^\infty b(y)dy\right) a(x)\\
&\leq-\bigg[M_0(\theta C_0+\theta\eps\|b'(x)\|_{L^{\infty}}t-\frac{1-\theta}{2}\|b\|_{L^1})\\
&\qquad-\frac{\theta}{4}(2C_0+2\eps\|b'(x)\|_{L^{\infty}}t+\|b\|_{L^1})\|b\|_{L^1}\bigg]a_0(x)\\
&\leq -M_0a_0(x)\bigg[\theta C_0-\frac{1}{2}(\theta C_0+(1-\theta)M_0+\frac{\theta}{2} M_0\|a_0\|_{L^1})\|a_0\|_{L^1}\bigg]\\
&\leq0.
\end{split}
\end{equation*}
The last inequality holds on the condition that
\begin{equation}\label{1}
\|a_0(x)\|_{L^1}\leq\frac{2\theta C_0}{\theta C_0+M_0} \text{ and } \|a_0(x)\|_{L^1}\leq1.
 \end{equation}
Then we also have
\begin{equation*}
\begin{split}
R_2&\geq \eps b'(x)+\eps\|b'(x)\|_{L^{\infty}}\\
&\qquad +b(x)\left(-\frac{1-\theta}{2}\phi+\frac{1+\theta}{2}\psi\right)-\theta\frac{(\phi+\psi)(\phi-\psi)}{4}a(x)\\
&\geq b(x)\left(\theta C_0+\eps\theta\|b'(x)\|_{L^{\infty}}t-\frac{1-\theta}{2}\int_{-\infty}^xb(y)dy
+\frac{1+\theta}{2}\int^{\infty}_xb(y)dy\right)\\
&\qquad -\frac{\theta}{4}({2C_0+2\eps\|b'(x)\|_{L^{\infty}}}t+\|b\|_{L^1})
\left(\int_{-\infty}^xb(y)dy-\int_x^\infty b(y)dy\right) a(x)\\
&\geq\bigg[M_0(\theta C_0+\theta\eps\|b'(x)\|_{L^{\infty}}t-\frac{1-\theta}{2}\|b\|_{L^1})\\
&\qquad-\frac{\theta}{4}(2C_0+2\eps\|b'(x)\|_{L^{\infty}}t+\|b\|_{L^1})\|b\|_{L^1}\bigg]a_0(x)\\
&\geq M_0a_0(x)\bigg[\theta C_0-\frac{1}{2}(\theta C_0+(1-\theta)M_0+\frac{\theta}{2} M_0\|a_0\|_{L^1})\|a_0\|_{L^1}\bigg]\\
&\geq0.
\end{split}
\end{equation*}
Hence $\textbf{(C2)}$ is verified for $(\bar{w}, \bar{z})$. From \eqref{2} and \eqref{1}, $a_0$ must satisfy
\begin{equation}\label{e:a0}
\|a_0\|_{L^1}\leq\min\left\{1, \, \frac{2\theta C_0}{\theta C_0+M_0}, \, \frac{1-\theta}{\theta}-\frac{C_0+1}{M_0}\right\}.
\end{equation}
 Now we turn to choose $M_0$ and $C_0.$ Considering the initial values of approximate solutions, we shall choose $C_0$ large enough first such that
\begin{equation*}
C_0\geq\max\{\sup w(x, 0), -\inf z(x, 0)\},
\end{equation*}
and then we have $w(x, 0)\leq\phi(x, 0), z(x, 0)\geq-\psi(x, 0).$  One choice of $M_0$ is
$$M_0=C_0\frac{3\theta^2+\theta}{1-\theta},$$
and then
$$\frac{2\theta C_0}{\theta C_0+M_0}=\frac{1-\theta}{1+\theta} \text{ and }
\frac{1-\theta}{\theta}-\frac{C_0}{M_0}-\frac{1}{M_0}\geq\frac{1-\theta}{1+\theta}$$
if $M_0$ is large enough. Thus our condition \eqref{a01} on $a_0$ satisfies \eqref{e:a0}, which is \textit{the key reason for \eqref{a01}}.
Therefore,  an application of Lemma \ref{modified maximum} yields
$$\bar{w}(x, t)\le0,~~ \bar{z}(x, t)\ge0,$$
which  implies
\begin{equation*}
\begin{split}
&w(x, t)\leq\phi(x, t)\leq C_0+\|b\|_{L^1}+1= C,\\
&z(x,t)\geq-\psi(x,t)\geq-C_0-\|b\|_{L^1}-1= -C,
\end{split}
\end{equation*}
where we can see that $C$ is independent of time.   Hence we obtain
\begin{equation}\label{uniform}
0\le\rho(x, t)\leq C,~~ |m(x, t)|\leq C\rho(x, t).
\end{equation}

\textbf{Step 2. Lower bound of density.}  By \eqref{uniform}, we know that the velocity $u=\frac{m}{\rho}$ is uniformly bounded, i.e., $|u|\le C$. Then the lower bound of density can be derived by the method of \cite{Huang20171}. Set $v=\ln\rho$, and then we get a scalar equation for $v$
\begin{equation}\label{v}
v_t+v_xu+u_x=\eps v_{xx}+\eps v_x^2+a(x)u
\end{equation}
from which we have
\begin{equation*}
v=\int_{\R}G(x-y, t)v_0(y)dy+\int_0^t\int_{\R}(\eps v_x^2-v_xu-u_x+a(x)u)G(x-y, t-s)dyds,
\end{equation*}
where $G$ is the heat kernel satisfying
\[
\int_\R G(x-y, t)dy=1, ~~\int_\R|G_y(x-y, t)|dy\leq\frac{C}{\sqrt{\eps t}}.
\]
Then it follows that
\begin{equation*}
\begin{split}
v&=\int_{\R}G(x-y, t)v_0(y)dy+\int_0^t\int_{\R}(\eps v_y^2-v_yu-u_y+au)G(x-y, t-s)dyds\\
&\geq\int_{\R}G(x-y, t)v_0(y)dy+\int_0^t\int_{\R}uG_y(x-y, t-s)+(au-\tfrac{u^2}{4\eps})G(x-y, t-s)dyds\\
&\geq\ln\eps-\frac{Ct}{\eps}-\frac{C\sqrt{t}}{\sqrt{\eps}}:=-C(\eps, t).
\end{split}
\end{equation*}
Thus \begin{equation}\label{rholow}
\rho\geq e^{-C(\eps, t)}.
\end{equation}

From \eqref{uniform} and \eqref{rholow}, we get \eqref{bound}. The lower bound of density guarantees that there is no singularity in \eqref{isen-vis}. Then we can apply classical theory of quasilinear parabolic systems  to complete the proof of Theorem \ref{thm-isenthvis}.


\subsection{Convergence of approximate solutions}

In this section, we will provide the proof of Theorem \ref{mainisen}.  Since we are focusing on the uniform bound of $\rho$ and $m,$ in this section we assume $1<\gamma\leq2$ for simplicity. For the case $2<\gamma\leq3$, one can follow the  similar argument in \cite{Lions2} or \cite{Wang} to obtain  the same conclusions.

Denote $\Pi_T=\R\times [0, T]$ for any $T\in(0, \infty).$

\textbf{Step 1. $H^{-1}_{loc}$ compactness of the entropy pair.}
We consider
\begin{equation*}
\eta(\rho^\eps, m^\eps)_t+q(\rho^\eps, m^\eps)_x,
\end{equation*}
where $(\eta,q)$ is  any weak entropy-entropy flux pair given in \eqref{2.6}. We will apply the Murat lemma to achieve the goal.
 \begin{lemma}{(Murat \cite{Murat})}\label{murat}
	Let $\Omega\in\R^n$ be an open set, then
	\begin{equation*}
		(\text{compact set of } W^{-1, q}_{loc}(\Omega))\cap(\text{bounded set of } W^{-1, r}_{loc}(\Omega))\\
		\subset(\text{compact set of } H^{-1}_{loc}(\Omega)),
	\end{equation*}
	where $1<q\leq 2<r.$
\end{lemma}

Let $K\subset\Pi_T$ be any compact set, and choose
$\varphi\in C_c^\infty(\Pi_T)$ such that $\varphi|_{K}=1$ and $0\leq\varphi\leq1.$
Multiplying \eqref{isen-vis} by  $\varphi\nabla\eta^* $ with $\eta^*$ the mechanical entropy,
we obtain
\begin{equation}\label{4.1}
\begin{split}
&\eps\int\int_{\Pi_T}\varphi(\rho_x, m_x)\nabla^2\eta^*(\rho_x, m_x)^\top dxdt\\
=&\int\int_{\Pi_T}(a(x)\frac{m^2}{\rho}-2\eps b(x)\rho_x)\eta^*_m\varphi+a(x)m\eta^*_\rho\varphi
+\eta^*\varphi_t+q^*\varphi_x+\eps\eta^*\varphi_{xx}dxdt.
\end{split}
\end{equation}
A direct calculation tells us that
\[
(\rho_x, m_x)\nabla^2\eta^*(\rho_x, m_x)^\top=p_0\gamma\rho^{\gamma-2}\rho_x^2
+\rho u_x^2.
\]
Noting that
\begin{equation*}
|(a(x)\frac{m^2}{\rho}-2\eps b(x)\rho_x)\eta^*_m|\leq\frac{\eps p_0\gamma}{2}\rho^{\gamma-2}\rho_x^2+\eps Cb^2m^2\rho^{-\gamma}+a_0\frac{m^{3}}{\rho^{2}},
\end{equation*}
we get
\begin{equation*}
\begin{split}
&\frac{\eps}{2}\int\int_{\Pi_T}\varphi(\rho_x, m_x)\nabla^2\eta^*(\rho_x, m_x)^\top dxdt\\
&\leq\int\int_{\Pi_T}(C\eps b^2m^2\rho^{-\gamma}+a_0\frac{m^{3}}{\rho^{2}})\varphi+
\eta^*\varphi_t+q^*\varphi_x+\eps\eta^*\varphi_{xx}\\
&~~+(\frac{m^{3}}{2\rho^{2}}
+\frac{\gamma}{\gamma-1}m\rho^{\gamma-1}p_0)a_0\varphi dxdt\\
&\leq C(\varphi).\\
\end{split}
\end{equation*}
Hence
\begin{equation}\label{en}\eps(\rho_x, m_x)\nabla^2\eta^*(\rho_x, m_x)^\top\in L^1_{loc}(\Pi_T),\end{equation}
i.e.,
\begin{equation}\label{isenlocestimate}
 \eps\rho^{\gamma-2}\rho_x^2+\eps\rho u_x^2\in L^1_{loc}(\Pi_T).
\end{equation}

For any weak entropy-entropy flux pairs given in \eqref{2.6}, as in \eqref{4.1}, we have
\begin{equation}\label{4.3}
\begin{split}
\eta_t+q_x&=\eps\eta_{xx}-\eps(\rho_x, m_x)\nabla^2\eta(\rho_x, m_x)^\top
+(\eta_\rho a(x)m+\eta_ma(x)\frac{m^2}{\rho})-2\eps \eta_m\rho_xb(x)\\
&=:\sum_{i=1}^4I_i.
\end{split}
\end{equation}
Using  \eqref{isenlocestimate},
it is straightforward to check that  $I_1$ is compact in $H^{-1}_{loc}(\Pi_T).$
Note that for any weak entropy, the Hessian matrix $\nabla^2\eta$ is controlled by  $\nabla^2\eta^*$ (\cite{Lions2}), that is,
\begin{equation}\label{4.4}
(\rho_x, m_x)\nabla^2\eta(\rho_x, m_x)^\top\leq(\rho_x, m_x)\nabla^2\eta^*(\rho_x, m_x)^\top,
\end{equation}
and thus $I_2$ is bounded in $L^1_{loc}(\Pi_T)$ and thus compact in
$W_{loc}^{-1, \alpha}(\Pi_T)$ for some $1<\alpha<2$ by the Sobolev embedding theorem.
For $I_3$, we have $$|I_3|=|\eta_\rho a(x)m+\eta_ma(x)\frac{m^2}{\rho}|\leq Ca_0,$$
which implies  that $I_3$ is bounded in $L^1_{loc}(\Pi_T).$
For the last term $I_4$, we get
$$|I_4|\leq C\eps \rho^{\gamma/2-1}|\rho_x|.$$ It follows from  \eqref{en} that $I_4$ is compact in $H^{-1}_{loc}(\Pi_T).$
Therefore,
$$\eta_t+q_x \text{ is compact in } W^{-1, \alpha}_{loc}(\Pi_T) \text{ with some } 1<\alpha<2.$$
On the other hand, since  $\rho$ and $m$ are uniformly bounded,  we have
$$\eta_t+q_x \text{ is bounded in } W^{-1, \infty}_{loc}(\Pi_T).$$
We conclude that
\begin{equation}\label{4.5}
\eta_t+q_x~ \text{is compact in}~ H^{-1}_{loc}(\Pi_T)
\end{equation}
for all weak entropy-entropy flux pairs with the help of the Murat lemma \ref{murat}.

\textbf{Step 2. Strong convergence and consistency.} By \eqref{4.5} and the compactness framework established in
\cite{Ding,DingDing, Diperna, Lions2},  we can prove that there exists
a subsequence of $(\rho^\eps,m^\eps)$ (still denoted by $(\rho^\eps,m^\eps)$) such  that
 \begin{equation}
 \label{4.6}
(\rho^\eps, m^\eps)\to(\rho, m) ~~~
  \text{ in } L^p_{loc}(\Pi_T), ~~p\geq1, \end{equation}
from which it is easy to show that $(\rho, m)$ is a weak  solution to the Cauchy problem \eqref{iso1}-\eqref{ini1}.  We omit the proof for brevity.

\textbf{Step 3. Entropy inequality.} We shall also prove that $(\rho, m)$ satisfies the entropy inequality in the sense of distributions for all weak convex entropies. Let $(\eta, q)$ be any entropy-entropy flux pair with $\eta$ being convex.  Multiplying \eqref{isen-vis} by $\varphi\nabla\eta $ with  $0\le \varphi\in C_c^\infty(\Pi_T)$,  we get
\begin{equation*}
\begin{split}
&\int\int_{\Pi_T}\eta_t\varphi+q_x\varphi dxdt\\
=&\int\int_{\Pi_T}\eps\eta_{xx}\varphi-\eps\varphi(\rho_x, m_x)\nabla^2\eta(\rho_x, m_x)^\top+\eta_{\rho}a(x)m\varphi+\eta_m(a(x)\frac{m^{2}}{\rho}-2\eps b(x)\rho_x)\varphi dxdt.
\end{split}
\end{equation*}
As in Step 1,  we have
\begin{equation*}
\left|\int\int_{\Pi_T}\eps\eta_{xx}\varphi dxdt\right|\to 0 \text{ as } \eps\to 0.
\end{equation*}
Moreover,
\begin{equation*}
\begin{split}
&\left|\int\int_{\Pi_T}2\eps \rho_xb(x)\eta_m\varphi dxdt\right|\\
\leq&\left[\int\int_{\Pi_T}C\eps \rho^{2-\gamma}b^2\varphi dxdt\right]^\frac{1}{2}
\left[\int\int_{\Pi_T}\varphi \eps\rho_x^2\rho^{\gamma-2}dxdt\right]^\frac{1}{2}\\
\leq&  C\eps^{\frac{1}{2}}\to 0 \text{ as } \eps\to0.
\end{split}
\end{equation*}
Noting that
$$\eps\varphi(\rho_x, m_x)\nabla^2\eta(\rho_x, m_x)^\top\geq0,$$
we conclude that
\begin{equation*}
\int\int_{\Pi_T}\eta\varphi_t+q\varphi_x dxdt
+(\eta_m\frac{m^{2}}{\rho}+m\eta_{\rho})a(x)\varphi dxdt\geq 0  \text{ as } \eps\to 0,
\end{equation*}
that is,  $(\rho,m)$ is indeed an entropy solution to the Cauchy problem \eqref{iso1}-\eqref{ini1}.
Therefore, the proof of Theorem \ref{mainisen} is completed.

\section{Preliminary and Formulation for Isothermal Flow}\label{formula2}
In this section, we provide some preliminaries and formulation for the isothermal case. Here, we adopt a similar notion as in Section \ref{formula} with no confusion. Letting
\begin{equation*}
n=A(x)\rho,~~ J=A(x)m,
\end{equation*}
and using $\gamma=1,$  we can rewrite \eqref{iso1} as
\begin{eqnarray}\label{iso2}
\left\{ \begin{array}{llll}
\displaystyle n_t+J_x=0,\\
\displaystyle J_t+\left(\frac{J^2}{n}+n\right)_x=-a(x)n,
x\in\R
\end{array}
\right.
\end{eqnarray}
with $a(x)=-\frac{A'(x)}{A(x)},$ $J=nu$. Then seeking weak entropy solutions of \eqref{iso1}-\eqref{ini1} is equivalent to solving \eqref{iso2} with the following initial data:
\begin{equation}\label{ini2}
(n, J)|_{t=0}=(n_0(x), J_0(x))=(A(x)\rho_0(x), A(x)m_0(x))\in L^{\infty}(\R).
\end{equation}
The eigenvalues of \eqref{iso2} are
\begin{equation*}
\lambda_1=\frac{J}{n}-1,\quad \lambda_2=\frac{J}{n}+1,
\end{equation*}
and the corresponding right eigenvectors are
\begin{equation*}
r_1=\left[\begin{array}{ll}
1\\ \lambda_1
\end{array}
\right],\quad
r_2=\left[\begin{array}{ll}
1\\ \lambda_2
\end{array}
\right].
\end{equation*}
The Riemann invariants $(w, z)$ are given by
\begin{equation*}
w=\frac{J}{n}+\ln n,\quad z=\frac{J}{n}-\ln n.
\end{equation*}
The mechanical energy $\eta^*(n, J)$ and mechanical energy flux $q^*(n, J)$ have the following formula
\begin{equation*}
\eta^*(n, J)=\frac{J^2}{2n}+n\ln n,~~ \quad q^*(n, J)=\frac{J^3}{2n^2}+J\ln n.
\end{equation*}

\section{Proof of Theorem \ref{main} }\label{isotheorem}

We first recall the compactness framework in Huang and Wang \cite{HuangWang}.

\begin{theorem}\label{framework}
Let $(n^\eps, J^\eps)$ be a sequence of bounded approximate solutions of \eqref{iso2}-\eqref{ini2} satisfying
\begin{equation*}
0<\delta\leq n^\eps\leq C, ~~|J^\eps|\leq n^\eps(C+|\ln n^\eps|)
\end{equation*}
with $C$ being independent of $\eps, T,$ $\delta=o(\eps)$. Assume that
\begin{equation*}
\partial_t\eta(n^\eps, J^\eps)+\partial_xq(n^\eps, J^\eps) \text{ is compact in } H^{-1}_{loc}(\Pi_T),
\end{equation*}
where $(\eta, q)$ is defined as
\begin{equation*}
\eta={n}^{\frac{1}{1-\xi^2}}e^{\frac{\xi}{1-\xi^2}\frac{J}{n}},~~
q=\left(\frac{J}{n}+\xi\right)\eta
\end{equation*}
for any fixed $\xi\in(-1, 1)$. Then there exists a subsequence of $(n^\eps, J^\eps)$, still denoted by $(n^\eps, J^\eps),$ such that
\begin{equation*}
(n^\eps(x, t), J^\eps(x, t))\rightarrow(n(x, t), J(x, t)) \text{ in } L^p_{loc}(\R\times\R^+),~ p\geq1,
\end{equation*}
for some function $(n(x, t), J(x, t))$ satisfying
\begin{equation*}
0\leq n\leq C, ~~|J|\leq n(C+|\ln n|),
\end{equation*}
where $C$ is a positive constant independent on $T.$
\end{theorem}

\subsection{Construction of approximate solutions}
Next we construct approximate solutions satisfying the conditions in Theorem \ref{framework}.  Raising density, which is motivated by \cite{lu1}, we add artificial viscosity as follows:
\begin{eqnarray}\label{iso-vis2}
\left\{ \begin{aligned}
\displaystyle &n_t+(J-\delta\frac{J}{n})_x=\eps n_{xx},\\
\displaystyle &J_t+\left(\frac{J^2}{n}-\frac{\delta}{2}\frac{J^2}{n^2}+\int^{n}_{\delta}\frac{t-\delta}{t}dt\right)_x=\eps J_{xx}-a(x)(n-\delta)+2b(x)\delta\frac{J}{n}-4\eps b(x)n_x
\end{aligned}
\right.
\end{eqnarray}
with initial data
\begin{equation}\label{ini-vis2}
(n, J)|_{t=0}=(n_0^\eps(x), J_0^\eps(x))=(n_0(x)+\delta, J_0(x))\ast j^\eps,
\end{equation}
where $b$ is a function to be determined later, $\delta=o(\eps),$ and $j^\eps$ is the standard mollifier and $0<\eps<1.$  By a direct computation, the eigenvalues are
\begin{equation}
\lambda^\delta_1=\frac{J}{n}-\frac{n-\delta}{n},\quad \lambda^\delta_2=\frac{J}{n}+\frac{n-\delta}{n},
\end{equation}
and the Riemann invariants are
\begin{equation*}
w=\frac{J}{n}+\ln n,\quad z=\frac{J}{n}-\ln n.
\end{equation*}
\subsection{Global existence of approximate solutions}\label{approximate-1}
%

In this section, we show the global existence of classical solutions to the Cauchy problem of quasilinear parabolic system \eqref{iso-vis2}-\eqref{ini-vis2} and obtain the following theorem.
\begin{theorem}\label{thm-isothvis}
There exists a unique global classical bounded solution $(n^\eps, J^\eps)$ to the Cauchy problem \eqref{iso-vis2}-\eqref{ini-vis2} satisfying \begin{equation}\label{isoth-bound}
\delta\leq n^\eps\leq C, ~~|J^\eps|\leq n^\eps(C+|\ln n^\eps|).
\end{equation}
\end{theorem}

We divide the proof of Theorem \ref{thm-isothvis} into three steps. In this section, we omit the up index $\eps.$

\textbf{Step 1. Local existence and lower bound of density.} The local existence of the solution for \eqref{iso-vis2}-\eqref{ini-vis2} can be proved by using the heat kernel and the same way in  \cite{Diperna}.
For the lower bound of density, we denote
$$v=n-\delta,$$
and then $v$ satisfies
\begin{equation}\label{v-1}
v_t+(uv)_x=\eps v_{xx},~~~~~~~~v|_{t=0}=v_0(x)
\end{equation}
with $u=\frac{J}{n}.$  From the definition of $n_0$, we have $v_0\geq0$.
Rewrite \eqref{v-1} as
\begin{align*}
v_t+uv_x=\eps v_{xx}-u_x v,
\end{align*}
and then it is easy to obtain from the maximum principle of the parabolic equation that
\begin{align*}
v(x, t)\geq\min v_0(x)e^{-\|u_x\|_{L^\infty}t}\geq0,
\end{align*}
and hence we gain $n\geq\delta.$

\textbf{Step 2. Uniform upper bound.} We apply Lemma \ref{modified maximum} to obtain the uniform $L^\infty$ estimates. As before, to estimate the uniform bound of the approximate solution, we shall investigate a parabolic system derived by Riemann invariants. We transform \eqref{iso-vis2} into the following form:
\begin{eqnarray*}
\left\{ \begin{aligned}
\displaystyle &w_t+\lambda^\delta_2 w_x=\eps w_{xx}+2\eps(w_x-2b(x))\frac{n_x}{n}-\eps\frac{n_x^2}{n^2}-a(x)\frac{n-\delta}{n}+2b(x)\delta\frac{J}{n^2},\\
\displaystyle &z_t+\lambda^\delta_1 z_x=\eps z_{xx}+2\eps(z_x-2b(x))\frac{n_x}{n}+\eps\frac{n_x^2}{n^2}-a(x)\frac{n-\delta}{n}+2b(x)\delta\frac{J}{n^2}.
\end{aligned}
\right.
\end{eqnarray*}
Set the control functions $(\phi,\psi)$ as follows:
\begin{equation*}
\begin{split}
&\phi=M+2\int_{-\infty}^xb(y)dy+2\eps\|b'(x)\|_{L^\infty}t,\\
&\psi=M+2\int^{\infty}_x b(y)dy+2\eps\|b'(x)\|_{L^\infty}t.
\end{split}
\end{equation*}
We remark that $\phi,\psi$ in this Section is different from those in Section \ref{isentheorem} for simplicity.
Then we obtain
\begin{equation*}
\begin{split}
&\phi_t=2\eps\|b'(x)\|_{L^\infty},~\phi_x=2b(x), ~\phi_{xx}=2b'(x);\\
&\psi_t=2\eps\|b'(x)\|_{L^\infty},~\psi_x=-2b(x), ~\psi_{xx}=-2b'(x).
\end{split}
\end{equation*}
Let $$\bar{w}=w-\phi,\bar{z}=z+\psi.$$
A simple calculation yields
\begin{eqnarray}\label{wz}
\left\{ \begin{aligned}
\displaystyle \bar{w}_t+\left(\lambda^\delta_2-2\eps\frac{n_x}{n}\right)\bar{w}_x=&\eps\bar{w}_{xx}+2\eps b'(x)-2\eps \|b'(x)\|_{L^\infty}-\eps\frac{n^2_x}{n^2}\\
&-2\left(\frac{J}{n}+\frac{n-\delta}{n}\right)b(x)-\frac{n-\delta}{n}a(x)+2b(x)\delta\frac{J}{n^2},\\
\displaystyle \bar{z}_t+\left(\lambda^\delta_{1}-2\eps\frac{n_x}{n}\right)\bar{z}_x=&\eps\bar{z}_{xx}+2\eps b'(x)+2\eps \|b'(x)\|_{L^\infty}+\eps\frac{n^2_x}{n^2}\\
&-2\left(\frac{J}{n}-\frac{n-\delta}{n}\right)b(x)-\frac{n-\delta}{n}a(x)+2b(x)\delta\frac{J}{n^2}.
\end{aligned}
\right.
\end{eqnarray}
Note that $$\frac{J}{n}=\frac{w+z}{2}=\frac{\bar{w}+\phi+\bar{z}-\psi}{2},$$
and then the system \eqref{wz} becomes
\begin{equation*}
\displaystyle\left\{ \begin{aligned} &\bar{w}_t+\left(\lambda^\delta_2-2\eps\frac{n_x}{n}\right)\bar{w}_x
=\eps\bar{w}_{xx}+a_{11}\bar{w}
+a_{12}\bar{z}+R_1,\\
&\bar{z}_t+\left(\lambda^\delta_{1}-2\eps\frac{n_x}{n}\right)\bar{z}_x
=\eps\bar{z}_{xx}+a_{21}\bar{w}
+a_{22}\bar{z}+R_2
\end{aligned}
\right.
\end{equation*}
with
\begin{equation*}
\begin{split}
&a_{11}=-b(x)\frac{n-\delta}{n},\quad ~~a_{12}=-b(x)\frac{n-\delta}{n}\leq0,\\
&a_{21}=-b(x)\frac{n-\delta}{n}\leq0,\quad ~~a_{22}=-b(x)\frac{n-\delta}{n}\\
\end{split}
\end{equation*}
and
\begin{equation*}
\begin{split}
R_1=&2\eps b'(x)-2\eps \|b'(x)\|_{L^\infty}-\eps\frac{n^2_x}{n^2}+(-a-b)\frac{n-\delta}{n}\\
&+b(x)\frac{n-\delta}{n}\left(2\int^{\infty}_xb(y)dy-2\int^x_{-\infty}b(y)dy-1\right),\\
R_2=&2\eps b'(x)+2\eps \|b'(x)\|_{L^\infty}+\eps\frac{n^2_x}{n^2}+(-a+b)\frac{n-\delta}{n}\\
&+b(x)\frac{n-\delta}{n}\left(2\int^{\infty}_xb(y)dy-2\int^x_{-\infty}b(y)dy+1\right),
\end{split}
\end{equation*}
where we have used $n\ge \delta$. Since $$\sup\left\{x\in\R\Big|2\int^{\infty}_xb(y)dy-2\int^x_{-\infty}b(y)dy\right\}=2\|b\|_{L^1},$$
we  can take $b(x)\in C^{1}(\R)$ such that
$$\|b(x)\|_{L^{1}}\leq\frac{1}{2}, \quad  |a(x)|\leq b(x), $$
and then we have $R_1\leq0,R_2\geq0.$ In fact, from our assumption on $a(x)$, we take $b(x)=a_0(x),$  which is \textit{our key reason for the condition \eqref{a02}}.   By our conditions on initial data,  we can take $M$ large enough such that
$$\bar{w}(x, 0)\le0,~~ \bar{z}(x,0)\ge0.$$
Then, Lemma \ref{modified maximum} yields
$$\bar{w}(x, t)\le0,~~ \bar{z}(x, t)\ge0,$$
which  implies that
\begin{equation*}
\begin{split}
&w(x, t)\leq\phi(x, t)\leq M+2\|b\|_{L^1}+2\eps \|b'\|_{L^{\infty}}t\leq C,\\
&z(x,t)\geq-\psi(x,t)\geq-M-2\|b\|_{L^1}-2\eps \|b'\|_{L^{\infty}}t\geq- C,
\end{split}
\end{equation*}
where for any fixed time $T,$ we choose $\eps$ small such that
$$\eps \|b'\|_{L^{\infty}}t\leq\eps \|b'\|_{L^{\infty}}T\leq1.$$
Hence we obtain \eqref{isoth-bound}.

From  Steps 1 and  2, using the classical theory of quasilinear parabolic systems, we can complete the proof of Theorem \ref{thm-isothvis}.
\subsection{Convergence of approximate solutions}\label{convergence}
%
 As stated in Section \ref{formula}, \eqref{iso1}-\eqref{ini1} is equivalent to \eqref{iso2}-\eqref{ini2}. Thus we only need to show that a subsequence of $(n^\eps, J^\eps)$ in Section \ref{approximate-1} converges to the solutions of \eqref{iso2}-\eqref{ini2} by verifying the conditions in Theorem \ref{framework}.  We also divide the proof into three steps.

\textbf{Step 1. $H^{-1}_{loc}$ compactness of the entropy pair.} We will verify the  $H^{-1}_{loc}$ compactness of the entropy pair
\begin{equation*}
\eta(n^\eps, J^\eps)_t+q(n^\eps, J^\eps)_x
\end{equation*}
for some weak entropy $(\eta, q)$ of \eqref{iso2} with
\begin{equation*}
\eta=n^{\frac{1}{1-\xi^2}}e^{\frac{\xi}{1-\xi^2}\frac{J}{n}},~~
q=\left(\frac{J}{n}+\xi\right)\eta
\end{equation*}
for any fixed $\xi\in(-1, 1).$
It is easy to calculate that
\begin{equation*}
\begin{split}
&\eta_n=\frac{1}{1-\xi^2}\left(1-\xi\frac{J}{n}\right)\frac{\eta}{n},~~\eta_J=\frac{\xi}{1-\xi^2}\frac{\eta}{n},\\
&\eta_{nn}=\frac{\xi^2}{(1-\xi^2)^2}\left(1-2\xi\frac{J}{n}+\frac{J^2}{n^2}\right)
n^{\frac{\xi^2}{1-\xi^2}-1}e^{\frac{\xi}{1-\xi^2}\frac{J}{n}},\\
&\eta_{nJ}=\frac{\xi^2}{(1-\xi^2)^2}\left(\xi-\frac{J}{n}\right)n^{\frac{\xi^2}{1-\xi^2}-1}e^{\frac{\xi}{1-\xi^2}\frac{J}{n}},\\
&\eta_{JJ}=\frac{\xi^2}{(1-\xi^2)^2}n^{\frac{\xi^2}{1-\xi^2}-1}e^{\frac{\xi}{1-\xi^2}\frac{J}{n}}.\\
\end{split}
\end{equation*}
Hence $$\eta_{nn}\eta_{JJ}-\eta_{nJ}^2=\frac{\xi^4}{(1-\xi^2)^3}n^{\frac{2\xi^2}{1-\xi^2}-2}e^{\frac{2\xi}{1-\xi^2}\frac{J}{n}}>0.$$
It indicates that $\eta$ is strictly convex for any $\xi\in(-1,1).$
Then
\begin{equation*}
\begin{split}
&(n_x, J_x)\nabla^2\eta(n_x, J_x)^\top\\
=&\frac{\xi^2}{(1-\xi^2)^2}n^{\frac{1}{1-\xi^2}-2}e^{\frac{\xi}{1-\xi^2}\frac{J}{n}}
\left[n_x^2+\left(\frac{J}{n}n_x-J_x\right)^2-2\xi n_x\left(\frac{J}{n}n_x-J_x\right)\right]\\
\geq&\frac{\xi^2}{(1-\xi^2)^2}\frac{\eta}{n^2}
\left[(1-|\xi|)n_x^2+(1-|\xi|)\left(\frac{J}{n}n_x-J_x\right)^2\right].
\end{split}
\end{equation*}
Let $K\subset\Pi_T$ be any compact set, and  choose
$\varphi\in C_c^\infty(\Pi_T)$ such that $\varphi|_{K}=1,$ and $0\leq\varphi\leq1.$
After multiplying \eqref{iso-vis2} by $ \varphi\nabla\eta,$
and integrating over $\Pi_T$, we obtain
\begin{equation*}
\begin{split}
&\eps\int\int_{\Pi_T}\varphi(n_x, J_x)\nabla^2\eta(n_x, J_x)^\top dxdt\\
=&\int\int_{\Pi_T}[-4\eps n_xb-a(n-\delta)+2b\delta\frac{J}{n}+\delta\frac{n_x}{n}+\frac{\delta}{2}(\frac{J^2}{n^2})_x]
\eta_J\varphi+\eta\varphi_t+\eps\eta\varphi_{xx}dxdt.\\
\end{split}
\end{equation*}
Due to
$$\left|\frac{J}{n}\right|\leq C+|\ln n|, \text{ and }\eta\leq n^{\frac{1}{1-\xi^2}}e^{\frac{|\xi|}{1-\xi^2}(C-\ln n)}\leq Cn^{\frac{1-|\xi|}{1-\xi^2}}, $$
it is easy to get
\begin{equation}\label{i3}
|-a(n-\delta)+2b\delta\frac{J}{n}\eta_J|\leq Cb.
\end{equation}
Besides,
\begin{equation}\label{i4}
|4\eps n_xb\eta_J|\leq\eps b|n_x|\frac{4|\xi|}{1-\xi^2}\frac{\eta}{n}\leq\frac{\eps\xi^2(1-|\xi|)}{4(1-\xi^2)^2}\frac{\eta}{n}\frac{n_x^2}{n}+C\eps\eta b^2.
\end{equation}
Moreover, we have
\begin{equation}
\begin{split}\label{i5}
&\left|\delta\frac{n_x}{n}\eta_J\right|\leq\frac{|n_x|}{n}\frac{|\xi|}{1-\xi^2}\frac{\eta}{n}\leq\frac{\eps\xi^2(1-|\xi|)}{4(1-\xi^2)^2}\frac{\eta}{n}\frac{n_x^2}{n}+C\frac{\delta^2}{\eps}\frac{\eta}{n^2},\\
&\left|\delta\left(\frac{J^2}{n^2}\right)_x\eta_J\right|\leq\left|\delta\frac{J}{n}\frac{|\xi|}{1-\xi^2}\frac{\eta}{n}\left(\frac{J}{n}\right)_x\right|\leq\frac{\eps\xi^2(1-|\xi|)}{4(1-\xi^2)^2}\eta\left|\left(\frac{J}{n}\right)_x\right|^2
+C\frac{\delta^2}{\eps}\frac{\eta}{n^2}\frac{J^2}{n^2}.
\end{split}
\end{equation}
Taking $\delta=\eps^3$ such that $\delta^2/\eps\leq\delta^{5/3}\leq n^{5/3},$ and choosing small $|\xi|\neq0$ , from the two facts
\begin{equation*}
\begin{split}
\frac{\eta}{n^2}&=n^{\frac{1}{1-\xi^2}-2}e^{\frac{\xi}{1-\xi^2}\frac{J}{n}}\leq Cn^{\frac{1}{1-\xi^2}-2-\frac{|\xi|}{1-\xi^2}}=Cn^{-\frac{2|\xi|+1}{1+|\xi|}},\\
\frac{\eta}{n^2}\frac{J^2}{n^2}&=n^{\frac{1}{1-\xi^2}-2}e^{\frac{\xi}{1-\xi^2}\frac{J}{n}}\frac{J^2}{n^2}\leq Cn^{\frac{1}{1-\xi^2}-2-\frac{|\xi|}{1-\xi^2}}(1+|\ln n|^2)\leq Cn^{-\frac{4|\xi|+1}{1+|\xi|}}
\end{split}
\end{equation*}
we get
\begin{equation}\label{i2}
\frac{\eps}{4}\int\int_{\Pi_T}\varphi(n_x, J_x)\nabla^2\eta(n_x, J_x)^\top dxdt\leq C(\varphi)
\end{equation}
with constant $C(\varphi)$ depending on the $H^2(\Pi_T)$ norm of $\varphi$. Hence for small $|\xi|\neq0,$
\begin{equation}\label{L1estimate}
\eps\frac{\eta}{n^2}n_x^2+\eps\frac{\eta}{n^2}\left(\frac{J}{n}n_x-J_x\right)^2
=\eps\frac{\eta}{n^2}n_x^2+\eps\eta\left|\left(\frac{J}{n}\right)_x\right|^2
\in L^1_{loc}(\Pi_T).
\end{equation}
Now we investigate the dissipation of the entropy as follows:
\begin{equation*}
\begin{split}
\eta_t+q_x=&\eps\eta_{xx}-\eps(n_x, J_x)\nabla^2\eta(n_x, J_x)^\top+[-a(n-\delta)+2b\delta\frac{J}{n}]\eta_J\\
&-4\eps n_xb\eta_J+\left(\delta(\frac{J}{n})_x\eta_n+[\delta\frac{n_x}{n}+\frac{\delta}{2}(\frac{J^2}{n^2})_x]\eta_J\right)\\
:=&\sum_{k=1}^5I_k.
\end{split}
\end{equation*}
Combining  \eqref{i3}, \eqref{i4}, \eqref{i5}, \eqref{i2}, we obtain that  $I_2+I_3+I_4+I_5$ is bounded in $L^1_{loc}(\Pi_T),$ and then compact in
 $W_{loc}^{-1, \alpha}(\Pi_T)$ with some $1<\alpha<2$ by the Sobolev embedding theorem.
For $I_1, $ from \eqref{L1estimate}, for any $\varphi\in H^1_0(\Pi_T),$
\begin{equation*}
\begin{split}
&\left|\int\int_{\Pi_T}\eps\eta_{xx}\varphi dxdt\right|
=\left|\int\int_{\Pi_T}\eps(\eta_nn_x+\eta_JJ_x)\varphi_x dxdt\right|\\
\leq&\int\int_{\Pi_T}\frac{\eps\eta|\varphi_x|}{1-\xi^2}
\left|\frac{n_x}{n}-\frac{\xi}{n}\left(\frac{J}{n}n_x+J_x\right)\right|dxdt\\
\leq&\sqrt{\eps}\left(\int\int_{\Pi_T}\frac{\eta\varphi_x^2}{n(1-\xi^2)}dxdt\right)^{\frac{1}{2}}
\bigg[\left(\int\int_{\Pi_T}\frac{\eps \eta n_x^2}{n^2}dxdt\right)^{\frac{1}{2}}\\
&~~~~~~+\left(\int\int_{\Pi_T}\frac{\eps\eta}{n^2}\left(\frac{J}{n}n_x-J_x\right)^2dxdt\right)^{\frac{1}{2}}
\bigg],
\end{split}
\end{equation*}
and thus we have that $I_1$ is compact in $H^{-1}_{loc}(\Pi_T).$
Finally, we get
$$\eta_t+q_x \text{ is compact in } W^{-1, \alpha}_{loc}(\Pi_T) \text{ with } 1<\alpha<2.$$
Moreover,
$$q=\left(\frac{J}{n}+\xi\right)\eta\leq (C-\ln n+|\xi|)\eta\leq
C+|\ln n|n^{\frac{1}{1-\xi^2}}e^{\frac{\xi}{1-\xi^2}\frac{J}{n}}\leq C,$$
and then
$$\eta_t+q_x \text{ bounded in } W^{-1, \infty}_{loc}(\Pi_T).$$
Therefore, taking $|\xi|$ small, we conclude that
$$\eta_t+q_x \text{ is compact in } H^{-1}_{loc}(\Pi_T) \text{ for small } |\xi|\leq 1,$$
by Lemma  \ref{murat}.

\textbf{Step 2. Convergence and consistency.} Since our  approximate solutions satisfy all the conditions in Theorem \ref{framework}, applying Theorem \ref{framework} yields
 $$(n^\eps, J^\eps)\to(n, J) ~~~
  \text{ in } L^p_{loc}(\Pi_T), ~~p\geq1.$$
This implies that $(n, J)$ is a weak solution to the Cauchy problem \eqref{iso2}-\eqref{ini2}. Similar to the previous argument, we can show that $(n, J)$ satisfies the energy inequality. Thus $(n, J)$ is an entropy solution. The proof of Theorem \ref{main} is completed.

\section{Appendix}
Here we provide the proof of Lemma \ref{modified maximum} for completeness.
\begin{proof}
Let $$\bar{M_0}=\|p\|_{L^\infty(\R\times[0,T])}+\|q\|_{L^\infty(\R\times[0,T])}.$$
We define two new variables
$$\bar{p}=p-\xi,~~
\bar{q}=q+\xi,$$
where $$\xi=\xi(x,t)=2\bar{M_0}\frac{\cosh x}{\cosh N}e^{\Lambda t}, N>0, $$
and $\Lambda>0$ will be determined later. For $\,(i,j)=(1,2)\,$or $\,(2,1),$ we write
\begin{equation*}
\begin{split}
a_{ij}(x,t,p,q)=&a_{ij}(x,t,\bar{p},\bar{q})\\
+&\left(\int_{0}^1 \frac{\partial a_{ij}}{\partial p}(x,t,\bar{p}+\tau\xi,\bar{q}-\tau\xi)d\tau-\int_{0}^1 \frac{\partial a_{ij}}{\partial q}(x,t,\bar{p}+\tau\xi,\bar{q}-\tau\xi)d\tau\right)\xi
\end{split}
\end{equation*}
and
\begin{equation*}
\begin{split}
\displaystyle &R_{i}(x,t,p,q,\zeta,\eta)=R_{i}(x,t,\bar{p},\bar{q},\zeta,\eta)\\
&+\left(\int_{0}^1 \frac{\partial R_{i}}{\partial p}(x,t,\bar{p}+\tau\xi,\bar{q}-\tau\xi,\zeta,\eta)d\tau-\int_{0}^1 \frac{\partial R_{i}}{\partial q}(x,t,\bar{p}+\tau\xi,\bar{q}-\tau\xi)d\tau\right)\xi.
\end{split}
\end{equation*}
Denote
\begin{align*}
&\overline{a}_{ij}=a_{ij}(x,t,\bar{p},\bar{q}),\overline{R}_{i}=R_{i}(x,t,\bar{p},\bar{q}),\\
&\overline{b}_{ij}=\int_{0}^1 \frac{\partial a_{ij}}{\partial p}(x,t,\bar{p}+\tau\xi,\bar{q}-\tau\xi)d\tau-\int_{0}^1 \frac{\partial a_{ij}}{\partial q}(x,t,\bar{p}+\tau\xi,\bar{q}-\tau\xi)d\tau,\\
&\overline{c}_{i}=\int_{0}^1 \frac{\partial R_{i}}{\partial p}(x,t,\bar{p}+\tau\xi,\bar{q}-\tau\xi,\zeta,\eta)d\tau-\int_{0}^1 \frac{\partial R_{i}}{\partial q}(x,t,\bar{p}+\tau\xi,\bar{q}-\tau\xi,\zeta,\eta)d\tau.
\end{align*}
Then we get a system for $(\bar{p}, \bar{q}),$
\begin{eqnarray*}
\left\{ \begin{aligned}
\displaystyle \bar{p}_t+\mu_1 \bar{p}_x=&
\eps \bar{p}_{xx}+a_{11}\bar{p}+\overline{a}_{12}\bar{q}+\overline{b}_{12}\bar{q}\xi+\bar{R}_1+\bar{c}_1\xi-2\mu_1\bar{M_0}\frac{\sinh x}{\cosh N}e^{\Lambda t}\\
&+\xi(-\Lambda+\eps+a_{11}-a_{12}),\\
\displaystyle \bar{q}_t+\mu_2 \bar{q}_x=&
\eps \bar{q}_{xx}+\overline{a}_{21}\bar{p}+a_{22}\bar{q}+\overline{b}_{21}\bar{p}\xi+\bar{R}_2+\bar{c}_2\xi+2\mu_2\bar{M_0}\frac{\sinh x}{\cosh N}e^{\Lambda t}\\
&+\xi(\Lambda-\eps+a_{21}-a_{22}).
\end{aligned}
\right.
\end{eqnarray*}
Note that  for any $\Lambda>0$,
$\bar{p}(x, t)< 0$ and $ \bar{q}(x, t)> 0$ hold for any $|x|\ge N$.
Next we show
\begin{align*}
\textbf{Claim}: &\text{There exists } \Lambda=\Lambda(\bar{M_0}) \text{ such that } \\
&\bar{p}(x, t)\le0 \text{ and } \bar{q}(x, t)\ge0 \text{ for } x\in(-N, N), 0\leq t\leq s_{\ast}=\frac{1}{\Lambda}.
\end{align*}
To this end, let
$$A=\{t\in[0,s_{\ast}]|\text{ there exist } x\in [-N,N] \text{ such that }\bar{p}(x, t)>0\text{ or }
\bar{q}(x, t)<0\}.$$
We shall prove the set $A$ is empty by contradiction. In fact, if $A$ is not empty, let $t_*=\inf A>0$, and then there exists $|x_*|\le N$ such that $\bar{p}(x_*, t_*)=0$ or $\bar{q}(x_*, t_*)=0.$
Without loss of generality, we assume $\bar{p}(x_*, t_*)=0.$
Then, $\bar{p}(x, 0)\leq0,~\bar{q}(x, 0)\geq0,\,|x|\leq N.$ For $0\leq t<t_*$,
$$~\bar{p}(\pm N, t)<0,~\bar{q}(\pm N, t)>0,~~|x|\le N,$$
and thus $\bar{p}(x, t)$ takes the maximum value over $[-N, N]\times[0, t_*]$
at the point $(x_*, t_*)$. We have
$$\bar{p}_x(x_*, t_*)=0,~\bar{p}_{xx}(x_*, t_*)\leq0,~\bar{p}_t(x_*, t_*)\geq0, ~\bar{q}(x_*, t_*)\geq 0.$$
Note that at the point $(x_*, t_*),$ $\overline{a}_{12}\leq0,\bar{R}_1\leq0, \Lambda t_{\ast}\leq \Lambda s_{\ast}=1.$
Moreover, for any $\tau\in[0,1],$ $$|\bar{p}+\tau\xi|\leq|p|+2\xi\leq\bar{M_0}+4\bar{M_0}e\leq C_1(\bar{M_0}),$$
$$|\bar{q}-\tau\xi|\leq|q|+2\xi\leq\bar{M_0}+4\bar{M_0}e\leq C_1(\bar{M_0}).$$
Therefore, $|\overline{b}_{12}|\leq C_2(\bar{M_0}),|\overline{b}_{21}|\leq C_2(\bar{M_0}),$ $|\overline{c}_{1}|\leq C_3(\bar{M_0}),|\overline{c}_{2}|\leq C_3(\bar{M_0}).$
A direct computation yields that at the point $(x_*, t_*)$,
\[ \begin{aligned}
\bar{p}_t+\mu_1\bar{p}_x \leq&\overline{a}_{12}\bar{q}+\overline{b}_{12}\bar{q}\xi+\bar{R}_1+\bar{c}_1\xi+\xi(-\Lambda+\eps+a_{11}-a_{12}+|\mu_1|)\\
\leq &\xi\left(-\Lambda+\eps+a_{11}-a_{12}+|\mu_1|+C_2(\bar{M_0})\bar{M_0}+2C_2(\bar{M_0})\bar{M_0}e+2C_3(\bar M_{0}) M_{0}e\right).\\
\end{aligned}\]
Then, choosing
$$\Lambda=:2\eps+\sum_{i=1}^{2}\|\mu_i\|_{L^\infty}+\sum_{i,j=1}^2 \|a_{ij}\|_{L^\infty}+C_2(\bar{M_0})\bar{M_0}(2e+1)+2C_3(\bar M_{0}) M_{0}e,$$
we get  $\bar{p}_t+\mu_1\bar{p}_x<0.$ It contradicts with
$$\bar{p}_t+\mu_1\bar{p}_x\geq 0,\text{ at }(x_{\ast},t_{\ast}).$$
Hence $A$ is empty and \textbf{Claim} holds. Letting
$N$ tend to infinity, we obtain that $ p(x,t)\leq 0,\, \mathrm{and}~q(x,t)\geq0 ~\, \forall x\in\R, 0\leq t\leq s_{\ast}.$
From the above analysis, we have proved that the set $$\Omega=\{t\in[0,T]|\,p(x,s)\leq0,\,\,q(x,s)\geq0 \,\forall x\in\R,\, 0\leq s\leq t\}$$ is an open set.
It is obvious that $\Omega$ is a closed subset of $[0,T].$ Therefore, $\Omega=[0,T]$. We thus complete  Lemma\,\ref{modified maximum}.
\end{proof}

\smallskip
\section*{Acknowledgments}
Wentao Cao's research is supported by ERC Grant Agreement No.724298. Feimin Huang is partially supported by National Center for Mathematics and Inter-disciplinary Sciences, AMSS, CAS, and NSFC Grant No.11371349 and 11688101. Difan Yuan is supported by China Scholarship Council No.201704910503. The authors would like to thank Professor Naoki Tsuge for valuable comments and suggestions.

\end{document}